\documentclass[11pt,reqno]{amsart}

\usepackage{latexsym}
\usepackage{amsfonts}
\usepackage{amsmath}
\usepackage{url}
\usepackage{amssymb,amsmath}

\def\bl{\begin{lemma}}
\def\el{\end{lemma}}
\def\bth{\begin{theorem}}
\def\eth{\end{theorem}}
\def\bc{\begin{corollary}}
\def\ec{\end{corollary}}
\def\bcj{\begin{conjecture}}
\def\ecj{\end{conjecture}}
\def\bpr{\begin{proposition}}
\def\epr{\end{proposition}}
\def\bde{\begin{definition}}
\def\ede{\end{definition}}

\def\H{\mathbb{H}}

\def\QED{\hfill\qedsymbol}

\newcommand{\be}{\begin{eqnarray}}
\newcommand{\ee}{\end{eqnarray}}

\newcommand{\R}{{\mathbb R}}

\newcommand{\N}{{\mathbb N}}

\renewcommand{\and}{\hbox{ {\rm and} }}

\newtheorem{theorem}{Theorem}[section]
\newtheorem{definition}{Definition}[section]
\newtheorem{lemma}[theorem]{Lemma}

\newtheorem{corollary}[theorem]{Corollary}
\newtheorem{proposition}[theorem]{Proposition}
\newtheorem{conjecture}[theorem]{Conjecture}
\newtheorem{quest}[theorem]{Question}

\theoremstyle{definition}
\numberwithin{equation}{section}

\begin{document}
\title{An Upper bound on the growth of Dirichlet tilings of hyperbolic spaces}
\author{Itai Benjamini \and Tsachik Gelander}
\date{}
\maketitle

\begin{abstract}
It is shown that the  growth rate $(\lim_r  |B(r)|^{1/r})$ of any $k$ faces Dirichlet tiling of $\H^d, d>2,$ is at most $k-1-\epsilon$,
for an $\epsilon > 0$,  depending only on $k$ and $d$. We don't know if there is a universal  $\epsilon_u > 0$, such that $k-1-\epsilon_u$ upperbounds the growth rate for  any $k$-regular tiling, when $ d > 2$? 
\end{abstract}

\section{Introduction}

Let $\H^d$ be the $d$-dimensional Lobachevsky space and let $G=\text{Isom}(\H^d)$ be its group of isometries.
Let $\Gamma\le G$ be a lattice and consider an associated {\em Dirichlet polyhedron} $\mathcal{P}$. Let $\mathcal{T}$ be the
tiling of $\H^d$ by $\Gamma$ translates of $\mathcal{P}$ and let $\mathcal{G}$ be the dual graph of $\mathcal{T}$.
That is, the vertices correspond  to the tiles and two vertices
are connected by an edge if they share a $d-1$ dimensional face.

Let
$$
gr(\mathcal{G}) = \lim_r  |B(r)|^{1/r},
$$
where $|B(r)|$ is the number of vertices in an $r$-ball of $\mathcal{G}$.
The limit exists since it is a submultiplicative sequence.

Suppose that $\mathcal{P}$ has $k$ faces, then the graph $\mathcal{G}$ is $k$-regular.\footnote{There are  only finitely many combinatorial types of polyhedra in dimension $d$ with $k$ faces, however each could a priory appear in infinitely many different tiling and the corresponding groups might be non-isomorphic.}

\bth
For any $d > 2$ and $k \in \N$, there is  a constant $\epsilon> 0$, such that $gr(\mathcal{G}) < k-1-\epsilon$ for any
$k$-regular graph dual to a finite volume convex tiling of $\H^d$.
\eth

Note that the $k$-regular tree, $T_k$, with $k\ge 4$ even, for which $gr(T_k) =  k-1$,
can be realized as a lattice graph in  $\H^2$.
\medskip

\begin{quest}
Is there a universal $\epsilon_u > 0$, so that any $k$-regular Dirichlet tiling $\mathcal{G}$ of  $\H^d, d >2$, $gr(\mathcal{G}) < k-1-\epsilon_u$?
\end{quest}

If such $\epsilon_u$ exists, possibly $k-1-\epsilon_u$ will upperbound the  growth rate of all tilings of high dimensional simply connected  homogenous spaces, which are not quasi $\H^2$. Also assuming it exists, then maybe the growth rate $k-1-\epsilon_u$ is attended for some tiling with a small $k$ in $\H^3$?  
\medskip

It is of interest to bound or calculate $\epsilon$ in the theorem. The argument below can be adapted to give some bounds.

{\it Remark:}
The set $S=\{\gamma\in \Gamma:\gamma\cdot\mathcal{P}~\text{is adjacent to}~\mathcal{P}\}$ generates $\Gamma$ and $\mathcal{G}$ is the corresponding Cayley graph
$\mathcal{G}=\text{Cay}(\Gamma,S)$.
\medskip

Let us note that in \cite{AGG} it is  shown that for any $\epsilon > 0$ there is a $4$-regular {\em amenable} Cayley graph
with  $gr(G)$ bigger than $3- \epsilon$.

\section{Proof}

\subsection{Dimension $d=3$}

\medskip


Here the idea is to show that the corresponding presentation
of the group will have some bounded relation, of length depending on $k$ only.
This implies that the growth rate  is uniformly smaller than for the free group.

Saying that there is a lattice whose fundamental domain is a
polytope with $k$ neighbours is the same as saying that the polytope has $k$ faces
and we have to specify how the faces are glued together. This means that
we have only a bounded number of ways $f(k)$ of gluing and
assuming it is  {\em torsion free}, then it follows that there is a bound on the
number of domains glued around an edge and the wanted bound follows.
\medskip

For the case of torsion, if we had an edge which has no
stabilizer then  around this edge there is a bound on the number of domains
touching it and hence a bound on a relation length and we are done. So
otherwise all edges have stabilizers which are just cyclic groups. But then consider
a vertex and two edges meeting at it.
We distinguish between the case where we could pick the vertex in the interior of the hyperbolic space and the case where all vertices are at infinity.

In the first case, we get in $SO(3)$ two cyclic
subgroups which generate a discrete group. These groups
are of bounded size as a finite subgroup of $SO(3)$ is either cyclic, dihedral or of order
at most $60$, (see e.g. \cite{M}). In our case the group cannot be dihedral since then one edge
stabilizer would be of order $2$.
If it is indeed so then again we are done.
\medskip

Thus we are left with the case of non-compact three dimensional lattices, where the fundamental
domain has a vertex at infinity. Here the corresponding group will
be inside a parabolic subgroup of $G$, hence solvable and by discreteness will preserve horospheres and hence must
be one of the finitely many non-cyclic discrete subgroups $\mbox{Isom}(\R^2)$,
and the argument continues like before.
\medskip

\subsection{Dimension $d\ge 4$}
Consider now the case of dimension strictly bigger than three.
By induction on $k$ one shows that there is a triangulation of bounded size
(easily computable) of the Dirichlet fundamental domain. Since hyperbolic simplices admit
an upper bound on their volume, this gives a bound on the volume. Now since $d >3$ it follows from
Wang's finiteness theorem \cite{W} that there are only finitely many possibilities for $\Gamma$. Thus we may suppose that $\Gamma$ is given.

 As in the discussion in dimension $3$ we may suppose that the stabilizer of every co-dimension $2$ face is non-trivial, and for each such face there is $\gamma$ in a bounded power $l$ of $S$ which stabilizes it.

Suppose first that $G/\Gamma$ is not compact.
In that case $\mathcal{P}$ has a vertex at infinity $\zeta$.
Let $\gamma_1,\gamma_2\in S^l$ be bounded elements stabilizing two distinct co-dim $2$ walls through $\zeta$ . Then
$\gamma_1,\gamma_2$ belongs to a spherical crystallographic subgroup of $\Gamma$,
hence by Bieberbach theorem (see e.g. \cite{V}) satisfies $[\gamma_1^m,\gamma_2^m]=1$ where $m$
depends only on $d-1$.
Thus, there is a relation of bounded length, which implies the growth gap in this case.

Next suppose that $G/\Gamma$ is compact. Here again one can argue similarly to the $3$ -dimensional case, however we give a different argument which can clearly be carried out in a much more general situation (when Wang's finiteness theorem applies).
Note that a Dirichlet domain is determined by a point,
and up to equivalence, it is enough to pick that point in a given fundamental domain, say $\mathcal{P_0}$.

Let $D=\text{diam}(\mathcal{P_0})$ and set $\tilde{\mathcal{P_0}}=\overline{N_D}(\mathcal{P_0})$, the closed $D$-neighborhood of $\mathcal{P_0}$.
Then $\tilde{\mathcal{P_0}}$ is compact, hence the set
$$
 \tilde{S}:=\{\gamma\in\Gamma:\gamma\cdot\tilde{\mathcal{P_0}}\cap\tilde{\mathcal{P_0}}\ne\emptyset\}
$$
is finite.

Any $\Gamma$-periodic $k$-regular graph $\mathcal{G}$ associated with a uniform Dirichlet  tiling of $\H^d$ is isomorphic to the Cayley graph $\text{Cay}(\Gamma:S)$ for some $S\subset\tilde{S}$, of cardinality
$k$. There are only finitely many such graphs, and these graphs are all transitive and not a tree. The result follows, since  $|B(r)|$ is a submultiplicative sequence.

\QED

\section{Further problems and remarks}

\begin{enumerate}

\item
Let  $L(k,d)$ be the set of $k$-regular dual graphs for Dirichlet tilings  in $\H^d$ and $gr(k,d)$ the set of growth rate  of these graphs.
Are there non-isomorphic such graphs with the same growth rate.
Bound the number of elements in these sets.  Are the values in $gr(k,d)$ rational or algebraic? How big/small can the exponents in $gr(k,d)$ be?
\medskip

\item
We conjecture that the maximal possible growth rate of $k$-regular dual graph  of a Dirichlet tiling,  is monotone decreasing in the dimension
(as long as $k > d$, for larger $k$  there is no such graph).

\item
A similar result should hold considering lattices in general real symmetric spaces of non-compact type of dimension greater than $2$, in particular the argument for compact Dirichlet domain that we gave for $d\ge 4$ applies in general.

\item
We suppose that a similar result should hold also for Cayley graphs, assuming non-amenability and vanishing of the first $L^2$ betti number, or property $T$.
What is the natural general statement?
\medskip

\item
In $\H^d$, $d>3$, are there only finitely many non-isomorphic such $k$-regular graphs for any $k$, estimate this number?
\medskip

\item
We studied here only Dirichlet tilings. We believe that similar bounds should hold for more general tilings. Including Voronoi tilings of point process
such as Poisson point process and aperiodic  tilings. See \cite{GS1, GS2, MS} for constructions of aperiodic tiling of the hyperbolic plane. There is no
construction yet in hyperbolic spaces of higher dimension and other symmetric spaces, such as $\H^2 \times \R$ or $\H^2 \times \H^2$.

\item
Given the tiling graph, from every vertex pick a geodesic to the root.
Do it by first connecting ball of radius $1$ then from all vertices at distance $2$ to vertices at distance $1$, and so on.
Look at the spanning subgraph consists of all vertices and only edges that are on such a chosen geodesic.
This graph has the same volume growth as the tiling graph.
Still maybe one can show directly that there is a positive portion of the edges that are not on this geodesic tree, at about every distance.
A direct argument might work  beyond Dirichlet tilings.

\end{enumerate}
\medskip

\noindent 
{\bf Acknowledgements:} we are grateful to Anton Malyshev  and Shahar Mozes.

\end{document}